\documentclass[a4paper,11pt]{article}
 \usepackage[english]{babel}
 \usepackage{graphicx}
 \usepackage[latin1]{inputenc}
 \usepackage[T1]{fontenc}
 \usepackage{lmodern}
 \usepackage[normalem]{ulem}
 \usepackage{verbatim}
 \usepackage{bbm}
 \usepackage{ntheorem}
 \usepackage{stmaryrd}
 \usepackage{amsmath}
 \usepackage{amssymb}
\usepackage{hyperref}

\theoremheaderfont{\scshape}
\theoremseparator{.}
\newtheorem{theorem}{Theorem}[section]

\newtheorem*{thmnonum}[theorem]{Theorem}

\def\sqw{\hbox{\rlap{\leavevmode\raise.3ex\hbox{$\sqcap$}}$%
\sqcup$}}

\newcommand{\N}{\ensuremath{{\mathbb N}}}

\newcommand{\Z}{\ensuremath{\mathbb Z}}

\newcommand{\Loc}{\textsc{l}}

\newcommand{\Gr}{\textsc{g}}

\newcommand{\Mon}{\textsc{i}}

\newcommand{\IH}{\textsc{ih}}

\newcommand{\Con}{\textsc{c}}

\newcommand{\pc}{{p_c}}

\newcommand{\SI}{\textsc{si}}

\newcommand{\defini}{\textbf}

\newcommand{\Conn}{\ensuremath{\mathsf{Conn}}}

\newcommand{\Cay}{\ensuremath{\mathsf{Cay}}}

\newcommand{\R}{\ensuremath{\mathbb R}}

\newcommand{\saut}{\vspace{0.25cm}}

\newenvironment{proof}{  
    \vspace*{-.4em}  {\it Proof.}%
}{
    \hfill\sqw\vspace*{.5em}
}

\author{S\'ebastien \sc{Martineau}\footnote{E-mail: sebastien.martineau@weizmann.ac.il.}\\
\\
{\it The Weizmann Institute of Science}}
\title{The set of connective constants of Cayley graphs contains a Cantor space}
\begin{document}
\maketitle

In this note, graphs are implicitly taken to be simple, unoriented, non-empty, connected and locally finite. Besides $\N$ is taken to contain 0.

\saut

A graph is said to be \defini{transitive} if it admits an action by graph automorphisms that is transitive on its set of vertices. Given a transitive graph $\mathcal{G}$ and a vertex $o$ of $\mathcal{G}$, denote by $c_n$ the number of paths starting at $o$, going through $n$ edges and not visiting any vertex more than once. By Fekete's Subadditive Lemma, the sequence $c_n^{1/n}$ converges to some real number $\mu(\mathcal{G})$. This number does not depend on the choice of $o$ and is called the \defini{connective constant} of $\mathcal{G}$.

\saut

Let us now define Cayley graphs. Given a group $G$ and a finite generating subset $S$ of $G$, the \defini{Cayley graph} associated with $(G,S)$ is the graph with vertex-set $G$ and such that two distinct elements $g$ and $h$ of $G$ are connected by an edge if and only if $g^{-1}h\in S\cup S^{-1}$. This defines a transitive graph $\Cay(G,S)$  which satisfies the implicit assumptions of this paper.

\saut

The purpose of this note is to prove the following theorem.

\begin{thmnonum}
The set $\{x:\exists (G,S),~x=\mu(\Cay(G,S))\}$ contains a Cantor space. In particular, this subset of $\R$ has cardinality $2^{\aleph_0}$.
\end{thmnonum}

This theorem implies the following result of Leader and Markstr\"om: the set of isomorphism classes of Cayley graphs has cardinality $2^{\aleph_0}$. See \cite{leadermark}.

\saut

An unpublished argument of Kozma \cite{gadykpersonal} shows that the set of $\pc$'s of Cayley graphs contains a Cantor space, where $\pc$ denotes the critical parameter for bond Bernoulli percolation. The strategy of proof used in the present paper is inspired by \cite{gadykpersonal}.

\vspace{1 cm}

\begin{proof}
Let $\Conn$ denote $\{x:\exists (G,S),~x=\mu(\Cay(G,S))\}$. Let $\Omega_\infty$ stand for $\{0,1\}^\N$, which is endowed with the product topology. It is enough to show that there is a continuous injection $f$ from $\Omega_\infty$ to $\Conn$. Indeed, as $\Omega_\infty$ is compact and $\R$ is Hausdorff, if there is such an $f$, then $f$ induces a homeomorphism from the Cantor space $\Omega_\infty$ to $f(\Omega_\infty)$. Besides, as $\Conn$ is a subset of $\R$ and both $\R$ and $\Omega_\infty$ have cardinality $2^{\aleph_0}$, the Cantor-Schr\"oder-Bernstein Theorem implies that if there is an injection from $\Omega_\infty$ to $\Conn$, then $\Conn$ has cardinality $2^{\aleph_0}$.

To prove the existence of a function $f$ as above, we will rely on several facts, which are listed below. Fact \Gr\ is about ``\underline{G}roups''. For it, the reader is referred to the \emph{proof} of Lemma III.40 in \cite{delaharpe}.
Facts \Loc\ and \SI, respectively on ``\underline{L}ocality'' and ``\underline{S}trict \underline{I}nequalities'', are due to Grimmett and Li: see respectively \cite{gllocality1} and \cite{glstrict}. Facts \Mon\ and \Con\ are easy and classical. They provide an ``\underline{I}nequality'' and a ``\underline{C}onvergence''.

\begin{itemize}
\item[\Gr] There are a finitely generated group $H$ and a subgroup $C$ of $H$ such that $C$ is isomorphic to $\bigoplus_{n\in \N} \Z$ and central in $H$.

\item[\Loc]  Let $(\mathcal{G}_n)_{n\leq \infty}$ be a sequence of Cayley graphs such that $\mathcal{G}_n$ converges locally\footnote{This means the following. Denote by $\rho_n$ the vertex corresponding to the identity element of $G_n$, where $\mathcal{G}_n=\Cay(G_n,S_n)$. For $n\leq \infty$ and $r\geq 0$, let $B_n(r)$ be the ball of centre $\rho_n$ and radius $r$ in $\mathcal{G}_n$, considered as a \emph{rooted graph}, rooted at $\rho_n$. Then we say that $\mathcal{G}_n$ \defini{converges locally} to $\mathcal{G}_\infty$ if $\forall r,~\exists n_0,~\forall n\geq n_0,~B_n(r)\simeq B_\infty(r)$.} to $\mathcal{G}_\infty$. Denote by $\mathcal{Z}$ the graph $\Cay(\Z,\{1\})$ and by $\mathcal{G}_n \times \mathcal{Z}$ the Cartesian product of the graphs $\mathcal{G}_n$ and $\mathcal{Z}$.

Then, $\mu(\mathcal{G}_n\times \mathcal{Z})$ converges to $\mu(\mathcal{G}_\infty \times \mathcal{Z})$.

\item[\SI] Let $G$ be a group generated by a finite subset $S$, and let $N$ be a normal subgroup of $G$. Assume that $N\not=\{1\}$ and that the ball of centre 1 and radius 2 of $\Cay(G,S)$ intersects $N$ only at 1.

Then, $\mu(\Cay(G/N,\overline{S})) < \mu(\Cay(G,S))$.

\item[\Mon] Let $G$ be a group generated by a finite subset $S$, and let $N$ be a normal subgroup of $G$.

Then, $\mu(\Cay(G/N,\overline{S})) \leq \mu(\Cay(G,S))$.

\item[\Con] Let $G$ be a group generated by a finite subset $S$. Let $(N_n)_{n\leq \infty}$ be a sequence of normal subgroups of $G$ such that for every finite subset $F$ of $G$, for $n$ large enough, $N_n \cap F= N_\infty \cap F$.

Then, $\Cay(G/N_n, \overline{S})$ converges locally to $\Cay(G/N_\infty, \overline{S})$.
\end{itemize}

\saut

The proof may now begin. Let us fix $(H,C_H)$ satisfying \Gr. Let $S_H$ be a finite generating subset of $H$. To avoid conflicts between additive and multiplicative notations, in the next three sentences, let us see $(\Z,+)$ as $(\langle a\rangle, \cdot)$. Let $G:= H \times \langle a\rangle$ and $S:= S_H \times \{1\} \cup \{(1,a)\}$. The finite subset $S$ of $G$ generates the group $G$. The subgroup $C:= C_H\times \{1\}$ of $G$ is central and isomorphic to $\bigoplus_{n\in \N} \Z$. Fix a basis  $(g_n)$ of the free abelian group $C$.

Let $\Omega$ denote the set of the (finite and infinite) words on the alphabet $\{0,1\}$. If $\mathcal{P}$ denotes a property which may be satisfied or not by elements of $\N\cup \{\infty\}$, denote by $\Omega_\mathcal{P}$ the set of the elements of $\Omega$ whose length satisfies $\mathcal{P}$. In this context, we may use the subscript ``$k$'' as an abbreviation for ``$=k$''. The $\Omega_\infty$ introduced at the beginning of the proof agrees with this notation.

For every $\omega\in \Omega_{<\infty}$, we will define a group $G_\omega$, which will be a quotient of $G$. Before stating our conditions, let us point out that we set $G_\text{empty word}$ to be $G$. As a result, $G$ with no subscript or with an empty subscript are both defined, and refer to the same object.

We will proceed by induction on $n\in \N$, with the following \underline{I}nduction \underline{H}ypothesis. See the figure on the next page.

\begin{itemize}
\item[\IH] For every $\omega \in \Omega_{\leq n}$, we have built a group $G_\omega$ which is $G$ or a quotient of $G$ by a subgroup of $\langle g_i:i<n\rangle$. We denote by $\mathcal{G}_\omega$ the Cayley graph of $G_\omega$ with respect to $\overline{S}$.

For every $\omega \in \Omega_{<n}$, we have constructed a real number denoted by $b_{\omega\star}$, and we have $G_{\omega0}=G_\omega$.

Setting $\text{``no letter''}<0<\star < 1$ and ordering lexicographically the words on the alphabet $\{0,\star,1\}$, the set
$$
\mathcal{S}_n := \{(\mu(\mathcal{G}_\omega),\omega):\omega\in \Omega_n\}\cup\{(b_{\omega\star},\omega\star):\omega\in \Omega_{<n}\}
$$
satisfies $\forall (x,\eta), (x',\eta')\in \mathcal{S}_n,~\eta < \eta' \iff x > x'$.
\end{itemize}
Notice that \IH\ holds for $n=0$.

\begin{figure}[h!]
\centering
\includegraphics{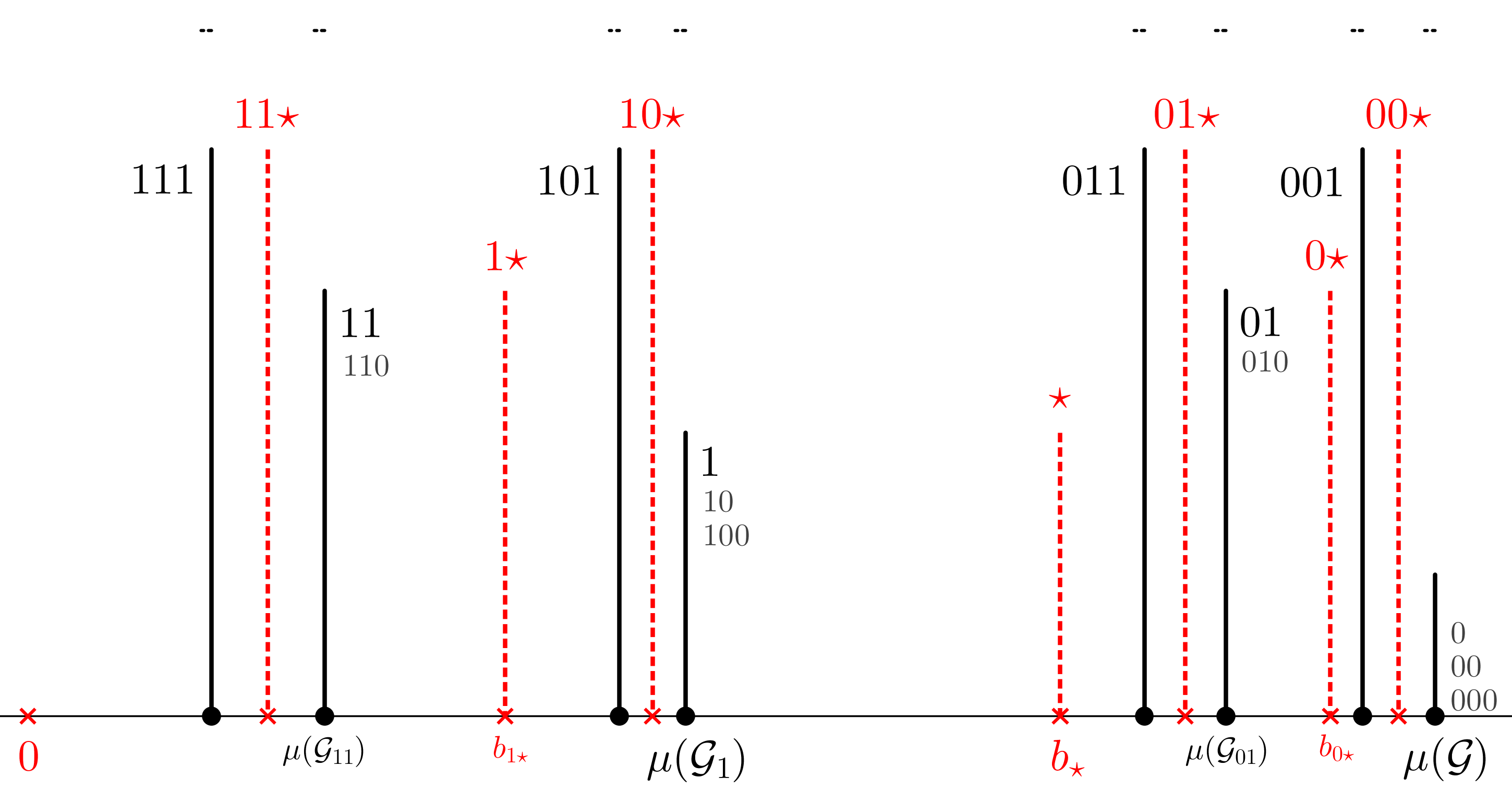}
{\textsc{Figure: }Illustration of \IH\ at rank 3. Above the vertical lines is ``represented'' the Cantor subset of $\Conn$ that we will build.}
\end{figure}

Let $n \in \N$ be such that \IH\ holds at rank $n$, and let us prove that it holds at rank $n+1$. 
For $\omega\in \Omega_n$ and $k\in \N^{\star}$, let $N_k^{\omega}$ denote the subgroup of $G_\omega$ generated by $\overline{g_n}^k$, which is normal in $G_\omega$ by centrality of $C$ in $G$.
Let $F$ be a finite subset of $G_\omega$. By \IH\ at rank $n$, the map $k\mapsto \overline{g_n}^k$ is injective from $\Z$ to $G_\omega$.
The set $Z_F^\omega:=\{j\in \Z: \overline{g_n}^j\in F\}$ is thus finite. For every $k > \max_{j\in Z_F^\omega} |j|$, we have $k\Z\cap Z_F^\omega \subset \{0\}$. As a result, for $k$ large enough $N_k^{\omega}\cap F = \{1\}\cap F$. It follows from \Con\ that $\Cay(G_\omega/\langle k \overline{g_n}\rangle,\overline{S})$ converges locally to $\mathcal{G}_\omega$ when $k$ goes to infinity.

Thus, by taking $m_n\in \N^\star$ large enough, \Loc\ and \SI\ guarantee that for every $\omega\in \Omega_n$, the connective constant $x:= \mu(\Cay(G_\omega/\langle m_n \overline{g_n}\rangle,\overline{S}))$ satisfies $x<\mu(\mathcal{G}_\omega)$ and, for every strict prefix $\alpha$ of $\omega$, $b_{\alpha\star}<x$. Taking $m_n$ to be minimal such that the above holds and letting $G_{\omega0}:= G_\omega$, $G_{\omega1}:=G_\omega/\langle m_n \overline{g_n}\rangle$ and $b_{\omega\star}:=(\mu(\mathcal{G}_{\omega0})+\mu(\mathcal{G}_{\omega1}))/2$, we get \IH\ at rank $n+1$. By induction, the $G_\omega$'s are constructed with the desired properties, together with the ``byproduct'' sequence $(m_n)$.

Now, for $\omega \in \Omega_\infty$, let $G_\omega := G/\langle m_i g_i : i \text{ such that }\omega(i)=1\rangle$ and let $\mathcal{G}_\omega:= \Cay(G_\omega, \overline{S})$.
To conclude the proof, it is enough to show that $f:\omega\mapsto \mu(\mathcal{G}_\omega)$ is injective and continuous as a function from $\Omega_\infty$ to $\R$.

Let $(\omega_n)$ be a converging sequence of elements of $\Omega_\infty$, and let $\omega_\infty$ denote its limit. For $n\in \N\cup\{\infty\}$, define $N_n$ to be $\langle m_i g_i : i \text{ such that }\omega_n(i)=1\rangle$. For every element $g=\prod_{i\leq i_0}g_i^{a_i}$ of the free abelian group $C$, we have
$$
g\in N_n \iff ``\{i:a_i\not=0\}\subset \{i:\omega_n(i)=1\}\text{ and }\forall i \leq i_0,~m_i | a_i\text{''}.
$$
Consequently, $(N_n)_{n\leq \infty}$ satisfies the hypotheses of \Con. By \Con\ and \Loc, $f(\omega_n)$ converges to $f(\omega_\infty)$, so that $f$ is continuous.

It remains to establish the injectivity of $f$. Let $\omega$ and $\omega'$ be two distinct elements of $\Omega_\infty$. Let $i\in \N$ be minimal such that $\omega(i)\not= \omega'(i)$. Without loss of generality, we may assume that $\omega(i)=0$ and $\omega'(i)=1$.
For $n\in \N$, denote by $\omega_n$ the prefix of $\omega$ of length $n$. Note that $\omega_i=\omega'_i$. By \Mon\ and the construction, we have
$$
\forall n>i,~\mu(\mathcal{G}_{\omega'})\leq \mu(\mathcal{G}_{\omega_i})< b_{\omega_i\star}< \mu(\mathcal{G}_{\omega_n}).
$$
By \Con\ and \Loc, $\mu(\mathcal{G}_{\omega_n})$ converges to $\mu(\mathcal{G}_\omega)$. Therefore, $$\mu(\mathcal{G}_{\omega'})\leq  \mu(\mathcal{G}_{\omega_i})< b_{\omega_i\star} \leq \mu(\mathcal{G}_{\omega}).$$ In particular, $\mu(\mathcal{G}_{\omega'}) \not= \mu(\mathcal{G}_{\omega})$. The function $f$ is thus injective, and the theorem is proved.
\end{proof}

\paragraph{Acknowledgements.} I would like to thank the Weizmann Institute of Science and my postdoctoral hosts --- Itai Benjamini and Gady Kozma --- for the excellent working conditions they have provided to me. I am also grateful to my postdoctoral hosts for letting me know of \cite{gadykpersonal}.

\begin{small}

\end{small}

\begin{thebibliography}{dlH00}

\saut

\bibitem[dlH00]{delaharpe}
Pierre de~la Harpe.
\newblock {\em Topics in geometric group theory}.
\newblock Chicago Lectures in Mathematics. University of Chicago Press,
  Chicago, IL, 2000.

\bibitem[GL]{gllocality1}
Geoffrey Grimmett and Zhongyang Li.
\newblock Locality of connective constants, {I}. {T}ransitive graphs.
\newblock ArXiv:1412.0150.

\bibitem[GL14]{glstrict}
Geoffrey Grimmett and Zhongyang Li.
\newblock Strict inequalities for connective constants of transitive graphs.
\newblock {\em SIAM Journal on Discrete Mathematics}, 28(3):1306--1333, 2014.

\bibitem[Koz]{gadykpersonal}
Gady Kozma.
\newblock Personal communication.

\bibitem[LM06]{leadermark}
Imre Leader and Klas Markstr{\"o}m.
\newblock Uncountable families of vertex-transitive graphs of finite degree.
\newblock {\em Discrete mathematics}, 306(7):678--679, 2006.

\end{thebibliography}
\end{document}